\documentstyle[12pt]{article}

\normalbaselines
\font\tbf = cmbx12
\pagebreak

\def\Clalgebra{$C_{p,N}$}

\begin{document}

\indent \vskip 0.6cm \centerline{\tbf EXTERIOR  DIFFERENTIALS OF
HIGHER  ORDER } \vskip 0.3cm \centerline{\tbf AND  THEIR
COVARIANT  GENERALIZATION} \vskip 1cm \centerline{\tbf V.
Abramov\footnote{\small Permanent address: Institute of Pure
Mathematics, University of Tartu, Vanemuise 46, Tartu, Estonia}
and R. Kerner} \vskip 1.5cm \centerline{ \it Laboratoire de
Gravitation et Cosmologie Relativistes} \vskip 0.2cm
\centerline{\it Tour 22, $4^{eme}$ \'etage, Bo\^{\i}te 142}
\vskip 0.1cm \centerline{\it Universit\'e Pierre et Marie Curie -
CNRS URA 769} \vskip 0.1cm \centerline{\it 4, Place Jussieu,
75005 Paris, France} \vskip 1.8cm \hskip 1cm {\tbf Abstract}
\vskip 0.2cm \indent {\small We investigate a particular
realization of generalized $q$-differential calculus of exterior
forms on a smooth manifold based on the assumption that $d^N=0$
while $d^k\neq 0$ for $k<N$. It implies the existence of cyclic
commutation relations for the differentials of first order and
their generalization for the differentials of higher order.
Special attention is paid to the cases $N=3$ and $N=4$. A
covariant basis of the algebra of such $q$-grade forms is
introduced, and the analogues of torsion and curvature of higher
order are considered. We also study a $Z_N$-graded exterior
calculus on a generalized Clifford algebra.}
\newpage
\noindent {\tbf I. INTRODUCTION} \vskip 0.3cm \indent An
appropriate framework for $d^N=0,\,N\geq 2$ generalization of
classical exterior differential calculus (satisfying $d^2=0$) is
provided by the notions of {\it graded $q$-differential algebra}
and {\it $q$-differential calculus} which have been elaborated in
recent papers ( \cite{Dubois-Kerner}, \cite{Dubois},
\cite{Majid}).

Let us remind that a graded $q$-differential algebra is a graded
unital algebra over the field $\bf C$ which is a sum of $N$
algebras with respective grade $k$, $k\in {\cal N}=\{0,1,2, \dots
N-1\}$: ${\cal A}=\oplus_{i\in {\cal N}}{\cal A}^i$, equipped
with an endomorphism $d$ of degree one satisfying the $q$-{\it
Leibniz rule}
$$d(AB)=d(A)\,B+q^a\,A\,d(B),\;A\in {\cal A}^a$$
and such that $d^N=0$ whenever $q^N=1$. A $q$-differential
calculus over an algebra ${\cal B}$ is a graded $q$-differential
algebra ${\cal A}$ such that $\cal B$ is a subalgebra of ${\cal
A}^0$.

These definitions show the way in which the exterior calculus of
differential forms on a smooth $n$-dimensional manifold $M$
should be generalized. The most striking property of this
generalized exterior calculus, due to the fact that $d^N=0$, is
that it contains not only first order differentials $dx^1, dx^2,
\dots, dx^n$ but must also include the {\it higher order
differentials} $d^2 x^1, d^2 x^1,\dots, d^{N-1}x^n$ as well.

After deriving from the condition $d^N f = 0$ a set of cyclic
commutation relations these differentials must satisfy, it
becomes clear that one needs {\it a generalization of the
Grassmann algebra} in addition to the above definitions in order
to produce a self-consistent algebra of generalized differential
forms. Such a generalization of Grassmann algebra which displays
a representation of the cyclic group $Z_3$ by cubic roots of
unity has been constructed in (\cite{Kerner1}) and then used in a
more general form in (\cite{Kerner2}, \cite{Kerner-Niemeyer},
\cite{AKR}) for the construction of the generalized exterior
calculus on a smooth manifold. It should be mentioned that
differential forms with higher order differentials have been
considered in \cite{Coquereaux}, where a formalism of
differential forms of higher order on any associative algebra has
been developed.

In this paper we continue to study a generalized $q$-exterior
calculus on a classical finite-dimensional smooth manifold $M$
paying particular attention to the tensorial behavior of the
generalized differential forms under a change of coordinates. The
main problem here is that the higher order differentials
transform in a non-homogeneous way under a coordinate
transformation. In order to circumvent this difficulty we
introduce an analogue of linear connection which allows us to
replace the ordinary differentials of any order by their
covariant generalizations.

The peculiar feature of linear connections we introduce is that
due to the higher order differentials its definition includes not
only the usual connection coefficients $\Gamma^i_{jk}$ but also
the additional coefficients $B^i_{jk},\;C^i_{jkl}$ (in the case
$N=3$) which {\it a priori} need not to be iterated by the first
order covariant differential. We find the transformation rules of
these coefficient functions of a connection which could be called
{\it connection coefficients of higher order} and we show that
$B^i_{jk}$ is related to the torsion of a connection. If there is
no torsion, then $C^i_{jkl}$ can be expressed in terms of the
Riemann curvature tensor.

We also study a particular realization of a $q$-exterior calculus
on a generalized Clifford algebra (\cite{Traubenberg},
\cite{Abramov}).
\newline
\vskip.4cm \noindent{\tbf II. GRADED ALGEBRA OF FIRST AND HIGHER
ORDER DIFFERENTIALS} \vskip 0.3cm \indent Let $M$ be a smooth
$n$-dimensional manifold and let $q$ be a $N$-th primitive root
of unity $q=e^{2\pi i/ N},\; q^N=1$. Let $U$ be an open subset of
$M$ with local coordinates $x^1, x^2, \ldots, x^n$. Our aim is to
construct an analogue of the exterior algebra of differential
forms with exterior differential $d$ satisfying the $q$-Leibniz
rule
\begin{equation}\label{Leibniz}
d(\omega \theta)=d\omega\,\theta + q^{\vert \omega\vert}\,\omega\,d\theta,
\end{equation}
where $\omega, \theta$ are complex valued differential forms,
$\vert \omega\vert$ is the degree of $\omega$, and
\begin{equation}\label{nilpotency1}
d^N=0,
\end{equation}
whereas $d^k\not =0$ for $1<k\leq N-1$. We shall also assume that
as in the classical case the exterior differential $d$ is a
linear operator and that it increases the degree of a form by one.

Let $U$ be an open subset of $M$ with the local coordinates
$x^{1},x^{2},\ldots ,x^{n}$. A differential form of degree zero is
a smooth function on $U$. Thus the set of differential forms of
degree zero $\Omega ^{0}(U)$ is the subalgebra of a whole algebra
which coincides with the algebra of smooth functions on $U$. A
differential 1-form on $M$ is an element of a free left module
$\Omega ^{1}(U)$ over the algebra $\Omega ^{0}(U)$ generated by
the differentials $dx^{1},dx^{2}, \ldots ,dx^{n}$ , and the right
module structure on $\Omega ^{1}(U)$ is defined by the relations
\begin{equation}
dx^{i}f(x)=f(x)\,dx^{i},\qquad f(x)\in \Omega ^{0}(U).
\label{module}
\end{equation}

The assumption $d^k\not =0$ for $1<k\leq N-1$ implies that there is no reason
to use only the first order differentials $dx^i$ in the construction of the
algebra of differential forms induced by $d$; one can also add a set of
formal {\it higher order differentials}, in which case the algebra will be
generated by
$$dx^1, \dots, dx^n, \dots, d^{N-1} x^1, d^{N-1} x^2, \dots, d^{N-1} x^n .$$
In order to endow the algebra of differential forms with appropriate
$Z_N$-grading we shall associate the degree $k$ to each differential
$d^k x^i$.  As usual, the grade of a product of differentials is the sum
of the degrees of its components modulo $N$. Given any smooth function $f$
and successively applying to it the exterior differential $d$ one obtains
the following expressions for the first three steps:
\begin{equation}\label{firstdifferential}
d f={\partial f\over \partial x^i}\; dx^i,\qquad
\end{equation}
\begin{equation}\label{seconddifferential}
d^2 f={\partial^2 f\over \partial x^i\partial x^j}\; dx^{(i} dx^{j)} +
{\partial f\over \partial x^i}\; d^2 x^i ,
\end{equation}
\begin{equation}\label{thirddifferential}
d^3 f={\partial^3 f\over \partial x^i\partial x^j\partial x^k}\;
dx^{(i} dx^j dx^{k)} + {\partial^2 f\over \partial x^i\partial x^j}\;
(d^2 x^i, dx^j)_q+ {\partial f\over \partial x^i}\; d^3 x^i.
\end{equation}
The relation (\ref{module}) between left and right structures of
the module of 1-forms $\Omega ^{1}(U)$ corresponds to classical
differential calculus on $U$ (\cite{Borowiec-Kharchenko}). Because
the partial derivatives of a smooth function of classical
differential calculus do commute, only the totally symmetric
combinations of indices are relevant in these definitions. That
is why in the above formulae the parentheses mean the
symmetrization with respect to the superscripts they contain, i.e.
\begin{eqnarray}
dx^{(i} dx^{j)}&=&{1\over 2!}\, (dx^i\,dx^j + dx^j\,dx^i),
dx^{(i}\nonumber\\
 dx^j dx^{k)}&=&\frac{1}{3!}\,(dx^i dx^j dx^k +
dx^j dx^k dx^i + dx^k dx^i dx^j\;\;\;\;\;\; \nonumber\\
        &&\;\;\;\;\;\;+ dx^k dx^j dx^i + dx^i dx^k dx^j + dx^j dx^i
        dx^k)\nonumber\\
(d^2 x^i, dx^j)_q&=&d^2 x^{(i}\, dx^{j)} +
         (1+q)\;d x^{(i} d^2 x^{j)}.\nonumber
\end{eqnarray}
The differentials of higher order of a function $f$ can be expressed by means
of a recurrent formula. Let us write the $k$-th differential of a function
$f$ in the form
\begin{eqnarray}\label{k-differential}
d^k f&=&\frac{\partial^{(k)} f}{\partial x^{i_1}\ldots \partial x^{i_k}}
L^{i_1 i_2\ldots i_k}_{(k)}+
 \frac{\partial^{(k-1)} f}{\partial x^{i_1}\ldots \partial x^{i_{k-1}}}
    L^{i_1 i_2\ldots i_{k-1}}_{(k)} \nonumber\\ &&\qquad + \ldots +
     \frac{\partial^{2} f}{\partial x^{i} \partial x^{j}}
        L^{ij}_{(k)}+\frac{\partial f}{\partial x^i} L^i_{(k)},
\end{eqnarray}
where $L^{i_1 i_2\ldots i_k}_{(k)}, \ldots, L^i_{(k)}$ are homogeneous
polynomials of differentials of total degree $k$ symmetric with respect to
their superscripts. They can be described by means of the following recurrent
formula
\begin{equation}\label{recurrent}
L^{i_1 i_2 \ldots i_m}_{(k)}=
        d L^{i_1 i_2 \ldots i_m}_{(k-1)}+\frac{1}{m}
             \sum^m_{l=1} dx^{i_l}
            L^{i_1 \ldots \hat{i_l}\ldots i_{m-1}}_{(k-1)},
\end{equation}
for $2\leq m\leq k-1$, and
\begin{equation}\label{aarmised}
L^{i_1 i_2\ldots i_k}_{(k)}=dx^{(i_1} dx^{i_2} \ldots dx^{i_k)},
\quad L^i_{(k)}=d^k x^i.
\end{equation}
In order to guarantee that the $N$-nilpotency (\ref{nilpotency1}) of the
$q$-exterior differential does not depend on the choice of local coordinates,
the $N$-th power of the differential $d$ should vanish identically on any
smooth function $f$ of a manifold $M$,
\begin{equation}\label{function}
d^N f=0.
\end{equation}
This leads to the conditions which should be imposed on formal differentials
$dx^1, dx^2, \ldots, \ldots, d^{N-1} x^1, d^{N-1} x^2, \ldots, d^{N-1} x^n$
in order to guarantee (\ref{function}). In their most general form they are
obtained from (\ref{k-differential}) and can be written as follows:
\begin{equation}\label{conditions}
L^{i_1 i_2 \ldots i_m}_{(k)}=0,\quad
    L^{i_1 i_2\ldots i_{k-1}}_{(k)}=0,\quad\ldots,\quad
         L^i_{(k)}=0.
\end{equation}
Let us write these conditions explicitly for the first few values of $N$.
If $N=2$ then (\ref{conditions}) takes on the form
\begin{equation}
L^{ij}_{(2)}=dx^{(i} dx^{j)}=0,\qquad L^i_{(2)}= d^2 x^i=0.
\end{equation}
Obviously these relations generate the classical exterior algebra based on
the skew-symmetric Grassmann structure with square nilpotent differential
$d^2=0$. The first non-trivial generalization of the classical algebra is the
case $N=3$ when the conditions (\ref{conditions}) take on the form :
\begin{eqnarray}\label{N=3}
L^{ijk}_{(3)}&=& dx^{(i} dx^j dx^{k)}=0,\cr
L^{ij}_{(3)}&=& d^2 x^{(i} dx^{j)}+(1+q) d x^{(i} d^2 x^{j)}=0,\cr
L^i_{(3)}&=& d^3 x^i = 0.
\end{eqnarray}
Although this paper concerns mainly with the $N=3$ generalization
of differential forms, we also show the constitutive relations
(\ref{conditions}) for $N=4$ :
\begin{eqnarray}\label{N=4}
L^{ijkl}_{(4)}&=& dx^{(i} dx^j dx^k dx^{l)}=0,\cr
L^{ijk}_{(4)}&=& d^2 x^{(i} dx^j dx^{k)}+
     (1+q)\; dx^{(i} d^2 x^j dx^{k)} +
         (1+q+q^2) dx^{(i} dx^j d^2 x^{k)} =0,\cr
L^{ij}_{(4)}&=&d^3 x^{(i} dx^{j)} +
     (1+q+q^2)\, d^2 x^{(i} d^2 x^{j)} +
            (1+q+q^2)\, d x^{(i} d^3 x^{j)},\cr
        L^i_{(4)}&=& d^4 x^i = 0.
\end{eqnarray}
The relations (\ref{conditions}) represent the minimal set of conditions that
should be imposed on the differentials in order to ensure (\ref{function}).
Comparing (\ref{conditions}) with (\ref{aarmised}) we conclude that for any
integer $N$ the differentials of first order $dx^i$ are $N$-nilpotent :
\begin{equation}\label{nilpotency2}
(dx^i)^N=0.
\end{equation}
On the other hand the relations (\ref{N=3}) and (\ref{N=4}) in special cases
of $N=3$ and $N=4$ demonstrate clearly that generally there are no relations
implying the nilpotency of any power for the differentials of higher order.
Therefore though the algebra generated by the relations (\ref{conditions}) is
finite-dimensional with respect to the first order differentials because of
(\ref{nilpotency2}), it remains infinite-dimensional with respect to the
entire set of differentials.

Since for $N>2$ the conditions (\ref{conditions}) do not
represent binary commutation relations, the algebra of
differential forms implemented by (\ref{conditions}) will be
rather hard to work with. One of the ways to circumvent this
difficulty is to find relations which on the one hand would be
simpler than (\ref{conditions}) but on the other hand they would
satisfy them. Following this idea we propose to solve the first
condition in (\ref{N=3}) by assuming that each cyclic permutation
of any three differentials of first order is accompanied by the
factor $q$ which in this case is a primitive cubic root of unity
and satisfies the identity
\begin{equation}
1+q+q^2=0.
\end{equation}
Thus we assume that each triple of differentials of first order
$dx^i, dx^j, dx^k$ is subjected to {\it ternary} commutation relations
\begin{equation}\label{ternary}
dx^i dx^j dx^k = q\; dx^j dx^k dx^i.
\end{equation}
These ternary commutation relations can not be made compatible with binary
commutation relations of any kind.

Therefore we suppose that all binary products $dx^i dx^j$ are
independent quantities. The second condition in (\ref{N=3}) can
be easily solved by assuming the following commutation relations:
\begin{equation}\label{binary}
dx^i d^2 x^l= q\; d^2 x^l dx^i.
\end{equation}
Note that from (\ref{ternary}) and (\ref{binary}) it follows that the above
ternary and binary commutation relations are coherent in the sense that
respect the grading defined earlier, i.e. the quantities $d x^k d x^m $ and
$d^2 x^j$ behave as elements of degree $2$ and could be interchanged in the
formulae (\ref{ternary}) and (\ref{binary}) .

The ternary commutation relations (\ref{ternary}) are much
stronger than the cubic nilpotence which follows from the first
relation of (\ref{N=3}). It has been proved in (\cite{Kerner2})
that if the generators of an associative algebra obey ternary
commutation relations such as (\ref{ternary}) then all the
expressions containing four generators should vanish. This means
that the highest degree monomials which can be made up of the
first order differentials have the form $dx^i dx^j dx^k\,, dx^i
(dx^j)^2$. Thus there are no fourth or higher degree differential
forms which can be made up of first differentials. In order to
construct an algebra with self-consistent structure we shall
extend this fact to the higher order differentials supposing that
{\it all} differential forms of fourth or higher degree vanish.

Since we have assumed that smooth functions commute with the
first order differentials (\ref{module}), i.e.
$$ x^k d x^m = d x^m x^k ,$$
then by virtue of the $q$-Leibniz rule
the second order differentials do not commute with smooth functions, because
differentiating the above equality we obtain
$$
d (x^k d x^m) = d x^k d x^m + x^k d^2 x^m = d ( d x^m x^k ) =
d^2 x^m x^k + q d x^m d x^k,
$$
which leads to the identity
\begin{equation}
x^k d^2 x^m - d^2 x^m x^k = q \, ( d x^k d x^m - q^2 \, dx^m d
x^k)
\end{equation}
In what follows, we shall consider only the expressions in which the forms
of different degrees are multiplied on the left by smooth functions of the
coordinates $x^k$, which means that we consider the algebra $\Omega (U)$
as a free finite-dimensional left module over the algebra of smooth
functions.

Let us find the number of independent generators $\cal N$ of this module. We
have $n$ first order differentials $dx^i$. The number of monomials spanning
the module of 2-forms is $n^2+n$ because we have $n^2$ independent binary
products $dx^i dx^j$ and $n$ second order differentials $d^2 x^i$. The number
of monomials spanning the module of 3-forms is $(n^3-n)/3 + n^2$ since there
are $(n^3-n)/3$ independent monomials $dx^i dx^j dx^k$ and $n^2$ independent
monomials $dx^i d^2x^j$. Summing all these numbers one finally obtains the
dimension of the module $\Omega (U)$
\begin{equation}
{\cal N} = \frac{n^3 + 6 n^2 + 5 n}{3}. \label{dimension}
\end{equation}
\indent
Although we have described the construction of the algebra $\Omega (U)$ only
in the case $N=3$ it can be extended to any integer $N>3$. In this case our
algebra is generated by the differentials
$dx^1,\ldots,dx^n,\ldots,d^{N-1}x^1,\ldots,d^{N-1}x^n$. Let
$d^{\alpha_1}x^{i_1} d^{\alpha_2}x^{i_2}\ldots d^{\alpha_r}x^{i_r}$ be a
monomial made up of differentials. We shall call the sum
$\alpha_1+\alpha_2+\ldots +\alpha_r$ an {\it order} of the monomial.
For the $N$-th order monomials we shall assume that they are subjected to
$r$-cyclic commutation relations
\begin{equation}\label{r-cyclic}
d^{\alpha_1}x^{i_1} d^{\alpha_2}x^{i_2}\ldots
d^{\alpha_{r-1}}x^{i_{r-1}}d^{\alpha_r}x^{i_r}=
 q^{\alpha_1}\;
            d^{\alpha_2}x^{i_2} d^{\alpha_3}x^{i_3}\ldots d^{\alpha_r}x^{i_r}
     d^{\alpha_1}x^{i_1},
\end{equation}
where $q$ is a $N$-th primitive root of unity. The relations
(\ref{ternary}) and (\ref{binary}), which determine the structure
of the algebra $\Omega(U)$ in the case of $N=3$, are the special
cases of the relations (\ref{r-cyclic}). We assume that the
monomials of order less than $N$ are independent. For the first
order differentials the $r$-cyclic relations (\ref{r-cyclic})
take on the form
\begin{equation}\label{N-cyclic}
dx^{i_1} dx^{i_2}\ldots
dx^{i_{N-1}}dx^{i_N}=
 q\;dx^{i_2} dx^{i_3}\ldots dx^{i_N} dx^{i_1}.
\end{equation}
Similarly to the case of $N=3$ it can be proved that the above $N$-cyclic
relations for first order differentials imply vanishing of all monomials
containing more than $N$ first order differentials. Extending this property
to the higher order differentials we shall assume that {\it all}
monomials of order higher than $N$ vanish.

In the next two Sections we show the examples of realization of
this exterior calculus. First we discuss the particular properties
of a $Z_3$-graded one- and two-dimensional realizations; then we
give an example of $p$ independent differentials acting on a
generalized Clifford algebra.
\newline
\vskip 0.4cm \noindent {\tbf III. EXAMPLES IN LOW DIMENSIONS}
\vskip 0.3cm \indent The aim of this Section is to investigate the
structure of the algebra of differential forms introduced above
by studying the simplest case of one-dimensional manifold and $N
= 3$. We shall denote the unique coordinate of this space by $t$.

Differentiating a smooth function $f$ one finds
\begin{eqnarray}
df&=&f'dt,\cr
d^2 f &=& f'' (dt)^2 + f' d^2 t,\cr
d^3 f &=& f''' (dt)^3 + f''(d^2 t dt+(1+q) dt d^2 t) +
             f' d^3 t.
\end{eqnarray}
In this simple case the above definitions yield immediately the relations
\begin{equation}
(dt)^3=0,\qquad dt\, d^2 t= q d^2 t\, dt.
\end{equation}
If one does not impose any additional relations, then the algebra of differential
forms based on the above commutation relations is infinite-dimensional and it
splits into the direct sum of two subspaces
\begin{eqnarray}
\Omega^{2m} = \{ \phi (d^2 t)^m + \psi (dt)^2 (d^2 t)^{m-1}\} ,  \, \ \ \, \ \
\Omega^{2m+1} = \{ \eta dt (d^2 t)^m \};
\end{eqnarray}
with  $ \phi, \, \ \ \psi \, \ \ \, \eta $ smooth functions of
$t$.

One has the following rules for calculating the exterior
differential
\begin{equation}
d (dt)=d^2 t,\qquad  d ( d^2 t)=d^3 t=0,
\end{equation}
\begin{equation}
\, {\rm and} \, \ \ d [(dt)^2] = d^2 t\, dt + q\;
dt\, d^2 t=(q+q^2) dt\, d^2 t=-dt\, d^2 t.
\end{equation}
\indent
It is interesting that the in their final form the rules for exterior
differentiation  do not contain the complex parameter $q$.

If $\omega\in \Omega^{2m}$ and
$\omega = \phi (d^2 t)^m + \psi (dt)^2 (d^2 t)^{m-1}$ then
$$ d\omega=(\phi' - \psi )\, dt (d^2 t)^m, $$
which means that $\omega$ is closed if and only if $\phi'=\psi$.
It is easy to show that any closed differential form of even degree
is exact. Indeed, if
$$\omega= \phi (d^2 t)^m + \phi' (dt)^2 (d^2 t)^{m-1} , $$
then $\omega=d\theta$, where $\theta\in \Omega^{2m-1}$ and
$$ \theta=\phi\;dt (d^2 t)^{m-1}. $$
From this it follows that for any differential form $\theta$ of odd degree
one has  $d^2 \theta = 0$.

Iterating twice the action of the exterior differential on an even degree
differential form $\omega= \phi (d^2 t)^m + \psi (dt)^2 (d^2 t)^{m-1}$ one
obtains the formula
$$ d^2 \omega= (\phi''-\psi') (dt)^2 (d^2 t)^m + (\phi'-\psi)(d^2 t)^{m+1}, $$
which shows that $d^2 \omega=0$ is equivalent to $d\omega=0$. Finally, if
$\theta=\eta\; dt (d^2 t)^m$ is an odd degree form, then
$$ d\theta=\eta' (dt)^2 (d^2 t)^m + \eta (d^2 t)^{m+1}, $$
and $d\theta=0$ implies $\eta'=\eta=0$. Thus any closed form of odd degree
is identically null. Now we turn to the transformation laws of differential
forms under the change of coordinates. Given a diffeomrphism $t=t (\tau)$
and a differential form of odd degree $\theta=\eta\; dt\,(d^2 t)^m$ one can
express it in coordinate $\tau$ as follows
$$ \theta=(t')^{m+1}\,\eta\;d\tau (d^2 \tau)^m, $$
which gives the transformation law for the coefficient function
$$ \eta=(t')^{m+1}\,\eta. $$
If $\omega=\phi (d^2 t)^m + \psi (dt)^2 (d^2 t)^{m-1}$ is a form of even
degree then expressing it in terms of coordinate $\tau$ one will obtain
$$ \omega=(t')^{m}\,\phi\;(d^2 \tau)^m + ([m]_q\,(t')^{m-1} t''\,\phi +
(t')^{m+1}\,\psi) (d\tau)^2 (d^2 \tau)^{m-1}, $$
where $[m]_q=1+q+q^2+\ldots + q^{m-1}$.
\newpage
\indent
The above formula gives the transformation law for the coefficient functions
of a form of even degree:
$$ \phi=(t')^m\,\phi,\qquad \psi=[m]_q\,(t')^{m-1} t''\phi +
(t')^{m+1} \,\psi. $$
We end this section by mentioning two facts. The first one is that given any
even differential form $\omega=\psi\,(dt)^2 (d^2 t)^{2l}$ one can solve the
equation $\theta^2=\omega$ by letting
$$\theta=q^{-l}\, \psi^{\frac{1}{2}}\,dt\,(d^2 t)^{l}.$$
\indent
We shall denote this solution by $\omega^{\frac{1}{2}}$. The second fact is
that given any $2m+1$- degree form $\theta=\eta\,dt\, (d^2 t)^m$ one can get
the closed $2m$- degree form by integrating with respect to $dt$, i.e. we
define the operator ${\cal I}_{ab}: \Omega^{2m+1}\to \Omega^{2m}$ by the
formula
$$ {\cal I}_{ab}(\theta)=(\int_a^b \eta\,dt)\,(d^2 t)^m, $$
where $a<b$ are finite real numbers. These facts gives us a possibility to
relate the differential forms we have described with the lenght of a smooth
curve on Riemannian manifold. Indeed let $M$ be a Riemannian manifold with
metric $g$ and  $\alpha: [a,b]\to M$ be a smooth curve on this manifold which
in local coordinates of the manifold $M$ is given by the equations
$x^i=x^i\,(t), \; a\leq t\leq b$. Then the first quadratic form
$ds^2=g_{ij}\,dx^i\,dx^j$ induces by means of the pullback the differential
2-form
$$ \omega=\alpha^{*}(ds^2)=g_{ij}\, \dot x^i\,\dot x^j\, (dt)^2. $$
If we now denote the length of a curve $\alpha$ by $S$ then
$$S={\cal I}_{ab}(\omega^{\frac{1}{2}}).$$
\indent
If we impose the vanishing of all monomials of degree $4$ and higher, then
on a two-dimensional real manifold with local coordinates $(x, y)$ the left
module of $Z_3$-graded forms has the dimension $14$, as it follows from the
general formula (\ref{dimension}). One has indeed to take into account the
following independent monomials:
$$ {\rm degree \ \ one : \ \ } \, \ \  dx, \ \ dy ; \ \ \, \ \ \,
{\rm degree \ \ two : \ \ } \, \ \ \, (dx)^2, \ \ (dy)^2, \ \  dx dy, \ \
dy dx,  \ \ d^2 x, \ \  d^2 y, $$
$$ {\rm and \ \ degree \ \ three : \ \ } \, \ \ \, d^2 x dx, \ \  d^2 x dy,
\ \ d^2 y dx, \ \  d^2 y dy, \ \  dx dx dy , \ \ dx dy dy $$
\indent
The particularity of the two-dimensional case is that it can be represented
in a more elegant way if we introduce complex notation with a single variable
$z = x + i \, y $. Then the $14$ independent {\it real} expressions above can
be expressed as
$$ d z = d x+i \, d y , \, \ \ d \bar{z} = \overline{d z } = d x-i \, d y, $$
$$ d^2 z = d^2 x + i \, d^2 y, \, \ \ d^2 \bar{z} = \overline{d^2 z} =
d^2 x - i \, d^2 y , $$
$$d z \, d z , \ \ d z \, d \bar{z} , \ \ d \bar{z} \,  d z =
\overline{ d z \, d \bar{z}}, \ \ d \bar{z} \, d \bar{z} =
\overline{ d z \, d z} $$
$$d^2 z \, d z, \ \ , d^2 \bar{z} \, d z, \ \ d^2 z \, d \bar{z} =
\overline{ d^2 \bar{z} \, d z}, \ \ d^2 \bar{z} \, d \bar{z} =
\overline{ d^2 z \, d z}$$
$$d z \, d z \, d z , \ \ \, {\rm and} \, \ \ d \bar{z} \, d \bar{z} \, d z
= \overline{ d z \, d z \, d \bar{z}} $$
\newpage
\indent
In two real dimensions the expression of a $1$-form $d f$, a $2$-form $d^2 f$
or of a $2$-form $\omega =d \theta$ with $\theta$ being an arbitrary $1$-form
are easily computed with the help of general formulae given in the Section $2$.

The situation becomes more interesting if we consider {\it complex holomorphic
functions} of the variable $z$. In such a case we have only one complex
variable and only two independent differentials, $d z$ and $d^2 z$ ; there is
no more need to introduce their complex conjugates. If we require now that
$d^3 f = 0$ for any holomorphic function $f$, then by virtue of
$$d f = \frac{d f}{d z} \, d z , \, \ \ d^2 f = \frac{d^2 f}{d z^2} \, d z \,
d z + \frac{d f}{d z} \, d^2 z, \, \ \ \, {\rm and \ \ imposing} $$
\begin{equation}
d^3 f = \frac{d^3 f}{d z^3} \, d z \, d z \, d z +
\frac{d^2 f}{d z^2} \, [ \, d^2 z \, d z + j d z \, d^2 z + d z \, d z \, ]
+ \frac{d f}{d z} \, d^3 z = 0 ,
\end{equation}
we arrive at the conditions on the differentials $d z$ and $d^2 z$ which are
formally the same as the ones verified by the single real variable $t$:
\begin{equation}
d^3 z = 0, \, \ \ \, (d z)^3 = 0, \, \ \ \, d z \, d^2 z = j \, d^2 z \, d z
\end{equation}
It is easy to show that the above relations imply that similar ones are
verified by the real differentials  $d x$ and $d y$:
$$d x \, d^2 x = j \,  d^2 x d x, \, \ \ \, d x \, d^2 y = j d^2 y \, d x,
\, \ \ {\rm etc.} $$
As in the real case, the algebra of degree $1$ forms is finite, but there
is no reason to cut off the powers of the degree $2$ forms $d^2 z$.
We don't see however how the integration introduced for the real
one-dimensional case can be generalized to the complex plane.
\newline
\vskip 0.4cm
\noindent
{\tbf IV. EXAMPLE OF q-EXTERIOR CALCULUS ON GENERALIZED CLIFFORD
ALGEBRAS} \vskip 0.3cm\indent In this Section we shall briefly
describe the structure of the generalized Clifford algebra
(\cite{Traubenberg}, \cite{Abramov}) and construct a $q$-exterior
calculus with $p$ exterior differentials $d_1, d_2, \ldots, d_p$
each satisfying $d_k^N=0.$ {\it Generalized Clifford algebra} is
an associative algebra over the complex field whose generators
$\Gamma_1, \Gamma_2, \ldots, \Gamma_p$ obey the commutation
relations
\begin{equation}
\Gamma_i\,\Gamma_j = q_{ij}\;\Gamma_j\,\Gamma_i, \, \ \ \, \ \ {\rm with}
\, \ \ \, \ \  \Gamma_k^N = 1,
\end{equation}
where
$$ q_{ij}=\left\{\begin{array}{l} 1,\quad\;\;\; i=j\\ q,\quad\;\;\; i<j\\
q^{-1}, \quad i>j\\ \end{array}\right. $$
It can be proved that the above commutation relations imply the
generalized Clifford relation
\begin{equation}\label{generalized delta}
\{\Gamma_{i_1}, \Gamma_{i_2},\ldots, \Gamma_{i_N}\} =
N!\,\delta_{i_1 i_2\ldots i_N},
\end{equation}
where the braces $\{\;,\;\ldots\;,\;\}$ at the left hand side stand for the
sum of all permutations with respect to the subscripts $i_1,i_2,\ldots, i_N$
which we shall call the $N$-anticommutator and $\delta_{i_1 i_2\ldots i_N}$
is the generalized Kronecker symbol which equals 1 when all subscripts are
equal and 0 in all other cases. Let us denote the generalized Clifford algebra
by $C_{p,N}$. This algebra can be endowed with $Z_N$ -grading if as usual one
associates grade 1 to each generator $\Gamma_k$ and defines the grade of any
monomial as a sum of the grades of the generators it is composed of modulo
$N$. Then the generalized Clifford algebra splits into the direct sum
$$ C_{p,N}=\sum_{i=0}^{N-1} C^{(i)}_{p,N}, $$
where $C^{(i)}_{p,N}$ is a subspace of the elements of grade $i$. The
dimension of the vector space underlying the algebra $C_{p,N}$ is $N^p$. It
can be also proved that the generalized Clifford algebra with $p$ generators
is isomoprhic to the grade zero subalgebra of the generalized Clifford
algebra with $p+1$ generators, i.e.
$$ C_{p,N}\cong C_{p+1,N}^{(0)}. $$
The generalized Clifford algebras have a matrix representations which can be
described as follows. Let us introduce the $n\times n$  matrices
\begin{equation}
\sigma_1=\pmatrix{0 & 1 & 0 & \ldots & 0 \cr
            0 & 0 & 1 & \ldots & 0 \cr
                \vdots & \vdots & \ddots & \vdots \cr
            0 & 0 & 0 & \ldots & 1 \cr
            1 & 0 & 0 & \ldots & 0 \cr},
\sigma_3=\pmatrix{1 & 0 & 0 & \ldots & 0 \cr
            0 & q & 0 & \ldots & 0 \cr
                \vdots & \vdots & \ddots & \vdots \cr
            0 & 0 & 0 & \ldots & 0 \cr
            0 & 0 & 0 & \ldots & q^{N-1} \cr},
\end{equation}
and $\sigma_2=(\sqrt{q})\;\sigma_3\,\sigma_1$, where $\sqrt{q}$ is needed
only in the case when $N$ is an even integer. Let $k=E(p/2)$. Then the
generators of the generalized Clifford algebra $C_{p,N}$ are represented by
the $n^k\times n^k$ matrices
\begin{eqnarray}\label{representation}
\Gamma_1 &=& \sigma_1\otimes I^{\otimes(k-1)},\qquad\qquad\qquad
            \Gamma_2=\sigma_2\otimes I^{\otimes(k-1)},\cr
            &\vdots&\qquad\qquad\qquad\qquad\qquad\qquad\qquad\vdots\cr
\Gamma_{2l-1} &=& \sigma_3^{\otimes(l-1)}\otimes\sigma_1
                     \otimes I^{\otimes(k-l-1)},\;\;\;
\Gamma_{2l}=\sigma_3^{\otimes(l-1)}\otimes\sigma_2
                     \otimes I^{\otimes(k-l-1)},\cr
         &\vdots&\qquad\qquad\qquad\qquad\qquad\qquad\qquad\vdots\cr
     \Gamma_{2k-1} &=& \sigma_3^{\otimes(k-1)}\otimes\sigma_1,\qquad
     \qquad\qquad
\Gamma_{2k}=\sigma_3^{\otimes(k-1)}\otimes\sigma_2,\cr
 & &\qquad\qquad\qquad \Gamma_{2k+1}=\sigma_3^{\otimes k},
\end{eqnarray}
where $I$ is the unit $N\times N$-matrix.

Because the generalized Clifford algebra $C_{p,N}$ possesses a natural
$Z_N$-grading one can use the $q$-commutator which is defined by the formula
$$ [B,B']_q=B\,B' - q^{bb'}\;B'\,B, $$
where $B,B'\in C_{p,N}$ and $b,b'$ are the grades of $B,B'$. Then
$q$-{\it exterior differentials} $d_1, d_2,\ldots, d_p$ are defined by the
formulae
\begin{equation}\label{q-differentials}
d_k B = [\Gamma_k, B]_q \equiv \Gamma_k\,B - q^{b} B\,\Gamma_k.
\end{equation}
According to the definition of $Z_n$-grading the $q$-exterior differential
raises the degree of an element $B$ by 1, i.e. $d_k:
C^{(i)}_{p,N}\to C^{(i+1)}_{p,N}$. It can be proved that each $q$-exterior
differential defined in (\ref{q-differentials}) is $N$-nilpotent
$$ d_k^N=0,\quad \mbox{for}\;k=1,\ldots, p. $$
Indeed if one writes the $l$-th power of the $q$-exterior differential
$d_k$ in the form
$$ d^l_k B=\sum_{i=0}^l \alpha^{(l)}_i \Gamma^{l-i}_k\,B\, \Gamma^i_k, $$
then the coefficients $\alpha^{(l)}_i$ are found to be
$$ \alpha^{(l)}_i=(-1)^i\,q^{\sigma}\,[l-i+1]^{(i-1)}_q,\qquad
         \sigma=\frac{2a+i-1}{2}\,i, $$
and
\begin{eqnarray}
[l]_q &=& 1+q+q^2+\ldots +q^{l-1},\cr
[l]'_q &=& 1+q\,[2]_q+q^2\,[3]_q+\ldots +q^{l-1}\,[l-1]_q,\cr
[l]''_q &=& 1+q\,[2]'_q+q^2\,[3]'_q+\ldots +q^{l-1}\,[l-1]'_q,\cr
&&\qquad\qquad\ldots\ldots\cr
[l]^{(i)}_q &=& 1+q\,[2]^{(i-1)}_q+q^2\,[3]^{(i-1)}_q+\ldots
                             +q^{l-1}\,[l-1]^{(i-1)}_q.
\end{eqnarray}
Thus in order to prove $N$-nilpotence of $q$-exterior differentials
suffice it to show that the relation
$$ [N-i+1]^{(i-1)}_q=0, $$
holds for every $i$ from 1 to $N-1$. But this is very easily proved by the
mathematical induction with respect to $i$.

It can be also proved that
\begin{equation}\label{differentials}
\{d_{i_1}, d_{i_2}, \ldots, d_{i_N}\}=0,\quad \mbox{for}\;\;
        1\leq i_1\leq {i_2}\leq \ldots\leq {i_N}\leq p.
\end{equation}
The above relations follow from (\ref{generalized delta}).

{\it The covariant differentials} $D_1, D_2, \ldots, D_p$ can be
defined by means of $q$-exterior differentials as follows
$$ D_k B=d_k B +A_k\,B, $$
where $A_k$ is a degree 1 element of the generalized Clifford algebra which
we shall call a $k$-{\it component of a connection 1-form} and use the
notation $A=(A_1, A_2, \ldots, A_p)$ combining all components into the
connection $A$.

Now if we apply the operator $\{D_{i_1}, D_{i_2}, \ldots, D_{i_N}\}$ to an
arbitrary element $B$ of the algebra the relations (\ref{differentials})
suggest that we get $B$ multiplied by an element of grade zero of the
algebra \Clalgebra which we call a $(i_1,i-2,\ldots,i_N)$-{\it component
of a curvature} and denote by $\Omega_{i_1 i_2\ldots i_N}$,
i.e.
\begin{equation}\label{curvature}
\{D_{i_1}, D_{i_2}, \ldots, D_{i_N}\}(B)=\Omega_{i_1 i_2\ldots i_N}\;B.
\end{equation}
Before giving the explicit expression for $\Omega_{i_1 i_2\ldots i_N}$ in
terms of connection in a general case we show the expressions for components
of curvature in low-dimensional cases of $N=2, 3$ and $p=2$. In the case of
$N=2, p=2$ the generalized Clifford algebra coincides with the classical
Clifford algebra represented by the Pauli matrices
$$ \sigma_1 = \pmatrix{0 & 1\cr 1 & 0\cr}, \quad  \sigma_2 =
\pmatrix{0& i\cr -i &0 \cr},\quad \sigma_3= \pmatrix{1&0 \cr 0& -1\cr}. $$
Computing the components of a curvature by means of the formula
(\ref{curvature}) one obtains
\begin{eqnarray}
\Omega_{11} &=& \{\sigma_1, A_1\}+A_1^2,\cr
\Omega_{12} &=& \{\sigma_1, A_2\}+\{\sigma_2, A_1\}+\{A_1,A_2\},\cr
\Omega_{22} &=&\{\sigma_2, A_2\}+A_2^2.
\end{eqnarray}
If the number of generators $p$ of the algebra remains the same but one takes
$N=3$ and denotes by $j$ a cubic root of unity to distinguish it from a
generic $N$-th root of unity $q$ then the components of curvature are
expressed in terms of the components of connection as follows
\begin{eqnarray}
\Omega_{111} &=& \{\eta_1, \eta_1, A_1\}+ \{\eta_1, A_1, A_1\}+
                A_1^3,\cr
\Omega_{112} &=& \{\eta_1, \eta_1, A_2\}+ \{\eta_1, A_1, A_2\}+
                \{A_1, \eta_1, \eta_2\}+ \{A_1, A_1, A_2\},\cr
\Omega_{122} &=& \{\eta_1, \eta_2, A_2\}+ \{\eta_1, A_2, A_2\}+
                \{A_1, \eta_2, \eta_2\}+ \{A_1, A_2, A_2\},\cr
\Omega_{222} &=& \{\eta_2, \eta_2, A_2\}+ \{\eta_2, A_2, A_2\}+
             +A^3_2,
\end{eqnarray}
where $\eta_1, \eta_2$ are the generators of the generalized Clifford
algebra $C_{2,3}$ and according to (\ref{representation}) they are
represented by the matrices
$$ \eta_1= \pmatrix{0 & 1 & 0\cr 0 & 0 & 1\cr 1 & 0 & 0\cr},\quad
 \eta_2= \pmatrix{0 & 1 & 0\cr 0 & 0 & j\cr j^2 & 0 & 0\cr}. $$
It is worth mentioning that the algebra generated by the above
matrices was dubbed by Sylvester the algebra of nonions
\cite{Sylvester}.

The above expressions for the components of a curvature in particular cases
$N=2, 3$ can be generalized for an arbitrary integers $p, N$ as follows. In
order to obtain the expression for the component $\Omega_{i_1 i_2 \ldots i_N}$
one should take the $N$-th anticommutator of generators
$\{\Gamma_{i_1}, \Gamma_{i_2}, \ldots, \Gamma_{i_N}\}$ and start replacing
the generators with the components of a connection with the same subscripts.
Let us introduce the following notations. Since there can be equal ones among
the integers $1\leq i_1\leq i_2\leq \ldots\leq i_N\leq p$ and they would give
us the same terms we pick only different ones denoting them by
$1\leq j_1\leq j_2\leq \ldots\ldots j_m\leq p$ and by $\vert j_k\vert$ the
number of integers in $(i_1, i_2, \ldots, i_N)$ equal to $j_k$. Then let us
denote by $\{A;\Gamma\}_{j_k}$ the anticommutator
$\{\Gamma_{i_1}, \Gamma_{i_2}, \ldots, \Gamma_{i_N}\}$ with $\Gamma_{j_k}$
being replaced with $A_{j_k}$, by $\{A;\Gamma\}_{j_k j_l}$ the same
anticommutator with $\Gamma_{j_k}, \Gamma_{j_l}$ being replaced with
$A_{j_k}, A_{j_l}$ and so on. It should be mentioned that subscripts $j_k$
and $j_l$ can be equal to each other if $\vert j_k\vert>1$. Then the
components of the curvature are expressed in terms of connection as follows:
$$\Omega_{i_1 \ldots i_N}=\{A;\Gamma\}_{j_1}+\ldots +\{A;\Gamma\}_{j_m}+
     \{A;\Gamma\}_{j_1 j_1} + \{A;\Gamma\}_{j_1 j_2}+
        \ldots + \{A;\Gamma\}_{j_1\ldots j_m} $$
\indent
The components of a curvature satisfy the Bianchi identities:
\begin{eqnarray}
 && d_{i_1}\Omega_{i_2 i_3\ldots i_{N+1}}+
 d_{i_2}\Omega_{i_1 i_3\ldots i_{N+1}}+
    \ldots + d_{i_{N+1}}\Omega_{i_1 i_2\ldots i_N}=\cr
        &&\qquad [\Omega_{i_2 i_3\ldots i_{N+1}}, A_{i_1}]_q+
        [\Omega_{i_1 i_3\ldots i_{N+1}}, A_{i_2}]_q+\ldots +
                [\Omega_{i_1 i_2\ldots i_N}, A_{i_{N+1}}]_q
\end{eqnarray}
\newline
\vskip 0.4cm \noindent {\tbf V. COVARIANT BASIS OF $Z_3$-GRADED
DIFFERENTIALS} \vskip 0.3cm \indent The $q$-exterior calculus on
generalized Clifford algebra described in the previous section
has a pure algebraic nature and though it is a good model of a
generalized exterior calculus with $d^N=0$ there even does not
arise the question of a change of coordinates. Going back to the
algebra $\Omega (U)$ introduced in the section 2 one might ask a
question whether this local algebra could be expanded on to the
whole manifold $M$. The difficulty here is that the set of
generators of the algebra includes the higher order differentials
which transform under a change of coordinates in a
non-homogeneous way. Our aim in this section is to show that
introducing an analogue of a linear connection and replacing the
ordinary differentials of all orders with the covariant ones we
can overcome the difficulty mentioned above (cf.
\cite{Kerner-Niemeyer}).

If we suppose that the formal expression $d^2 f$ does not vanish identically
as it is the case in the usual $Z_2$-graded exterior calculus of forms, then
we must abandon the antisymmetry of the product of $1$-forms in this algebra.
The vanishing of the expression (\ref{seconddifferential}) could be given
an intrinsic sense in any local coordinate system because one supposes that
simultaneously $d^2 = 0$, so it applied to any second differential of a local
coordinate, be it $d x^k$ or $d^2 y^{k'}$, and parallelly, taking into account
the symmetry of partial second derivatives,
$$\frac{\partial^2 f}{\partial x^k \partial x^m} =
\frac{\partial^2 f}{\partial x^m \partial x^k} ,$$
it had to be postulated that the product of $1$-forms must be antisymmetric:
$$ dx^k \, dx^m = - d x^m d x^k$$
Under a change of local coordinates, $x^k \rightarrow y^{m'} \, (x^k)$, and
$x^k = x^k( y^{m'}),$ the basis of $1$-forms transformed as a covariant tensor, i.e.
$d x^k \, = \, \frac{\partial x^k}{\partial y^{m'}} \, d y^{m'} . $
However, had we abandoned the postulate $d^2 = 0$ and the antisymmetry of the
product of $1$-forms, the {\it second differentials} of the coordinates, which
are for the time being purely formal expressions, would not transform as
covariant tensors because of the non-homogeneous term:
\begin{equation}\label{trans of sec diff}
d^2 x^k = d ( \frac{\partial x^k}{\partial y^{m'}} \, d y^{m'} ) =
\frac{\partial^2 x^k}{\partial y^{l'} \partial y^{m'}} \, d y^{l'} d y^{m'}
\, + \, \frac{\partial x^k}{\partial y^{m'}} \, d^2 y^{m'}
\end{equation}
Introducing {\it connection coefficients} $\Gamma^k_{lm}$ we shall define
the {\it covariant second differentials} $D^2 x^k$  as
\begin{equation}
D^2 x^k = d^2 x^k + \Gamma^k_{l m} \, d x^l \, d x^m
\end{equation}
(in order to make our notation homogeneous, from now on we shall also note
the first differentials - which are naturally covariant - with capital $D$,
i.e. we shall identify $D x^k = d x^k$.).
Note that the above equation can be still interpreted in terms of Grassmann
algebra of exterior forms: if we still impose $d^2 = 0$ and the antisymmetry
of the exterior product, the covariant $2$-form  $D^2 x^k$ is equal to the
torsion $2$-form. The vanishing of $D^2 x^k$ is then equivalent to the
condition of null torsion, which is valid in all coordinate systems.

Let us suppose now that the differentials $d x^k$ and $d^2 x^m$ satisfy the
relations imposed by the condition $d^3 = 0$ derived before, i.e., with
$q = e^{\frac{2 \pi i}{3}}$:
$$ d x^k \, d x^l \, d x^m = q \,  d x^l \, d x^m \, d x^k \, \ \ \,
\ \ {\rm and} \, \ \ \, \ \ d x^k \, d^2 x^m = q \,  d^2 x^m \, d x^k$$
It is obvious that these relations remain valid if we replace the ordinary
first and second differentials by their {\it covariant} counterparts:
\begin{equation}
D x^k \, D x^l \, D x^m = q \,  D x^l \, D x^m \, D x^k \, \ \ \,
\ \ {\rm and} \, \ \ \, \ \ D x^k \, D^2 x^m = q \,  D^2 x^m \, D x^k
\end{equation}
which span a covariant basis of the same $Z_3$-graded algebra,
which has the property of transforming covariantly under the
change of a basis.

Now, although $D^2 x^k$ represents a tensorial quantity, its
covariant differential $D \, (D^2 x^k) \, = \, D^3 \, x^k$ can
not be computed by simple iteration, i.e. by applying the same
formula as for the covariant differential of $D x^k$. As a matter
of fact, $D^3 x^k$ has to be a tensorial quantity of third
degree, which in covariant basis should contain both $D x^k \,
D^2 x^m$ and $D x^k D x^l D x^m$ . That is why we should write:
\begin{equation}\label{cube of D}
D^3 x^k = D \, (D^2 x^k) = d \, (D^2 x^k) + B^k_{lm} \, D x^l D^2 x^m +
C^k_{lmn} \, D x^l \, D x^m \, D x^n
\end{equation}
with new coefficients of two kinds, whose transformation properties under
coordinate change are yet to be derived, and which will assure the tensorial
character of $D^3 x^k$. Acting with the operator $d$ on $D^2 x^k$ yields the
explicit result:
\begin{eqnarray}
D^3 x^k &=& d^3 x^k + \partial_n \Gamma^k_{lm} \, d x^n d x^l d x^m +
\Gamma^k_{lm} \, d^2 x^l d x^m + q \, \Gamma^k_{lm} \, d x^l d^2 x^m +\cr
 && \qquad\quad+ B^k_{lm} \, D x^l D^2 x^m + C^k_{lmn} \, D x^l D x^m D x^n
\end{eqnarray}
$${\rm Now, \ \ using \ \ the \ \ fact \ \ that \ \ } \ \
d^2 x^l = D^2 x^l  - \Gamma^l_{rm} \, d x^r \, d x^m \, , $$
we can express $D^3 x^k$ by means of covariant quantities only:
\begin{eqnarray}
D^3 x^k &=& d^3 x^k + \biggl[ \, B^k_{lm} \, + \, q^2 \, \Gamma^k_{ml} + q \,
\Gamma^k_{lm} \, ] \, D x^l D^2 x^m  \, +\cr
&&\;\; +\biggl[ \, C^k_{lmn} \, + \, \partial_l \Gamma^k_{mn} - \Gamma^r_{lm}
\Gamma^k_{rn} \, - \, q \, \Gamma^r_{mn} \Gamma^k_{lr} \, \biggr] \, D x^l
D x^m D x^m \cr
&=& \, d^3 x^k + {\tilde{B}}^k_{lm} \, D x^l D^2 x^l +
{\tilde{C}}^k_{lmn} \, D x^l D x^m D x^n
\end{eqnarray}
Now we shall proceed by analogy with the $Z_2$-graded case. As we have seen,
the condition $d^3 = 0$ implies also the ternary and binary $q$-commutation
relations with $q = e^{\frac{2 \pi i}{3}}$. This means that in the final
expression we remain with
$$D^3 x^k = {\tilde{B}}^k_{lm} \, D x^l D^2 x^m + {\tilde{C}}^k_{lmn} \,
D x^l D x^m D x^n ,$$
It is obvious that if we want to impose the tensorial behavior on $D^3 x^k$,
then both coefficients ${\tilde{B}}^k_{lm}$ and ${\tilde{C}}^k_{lmn}$ must
transform as tensors given the manifestly tensorial character of the  products
of differentials they are contracted with. In contrast, the coefficients
$B^k_{lm}$ and $C^k_{lmn}$ have obviously non-tensorial character, which is
compensated by the connection coefficients and their derivatives entering
the definitions of ${\tilde{B}}^k_{lm}$ and ${\tilde{C}}^k_{lmn}$.
In order to get the transformation rules for the coefficients
$B^k_{lm}$ and $C^k_{lmn}$ we use the observation that the formula
(\ref{cube of D}) implicitly determines how the covariant differential $D$
is acting on the second order differentials. Indeed the left-hand side of
(\ref{cube of D}) can be written in the form
\begin{eqnarray}
D^3 x^k &=& D (D^2 x^k)=D ( d^2 x^k + \Gamma^k_{rs}\,Dx^r Dx^s)\cr
     &=& D ( d^2 x^k)+ \frac{\partial \Gamma^k_{rs}}{\partial x^l}\,Dx^l
             Dx^r Dx^s +
        (q\; \Gamma^k_{rs} + q^2\;\Gamma^k_{sr})\,dx^r D^2 x^s.
\end{eqnarray}
Before we express $D ( d^2 x^k)$ in terms of the coefficients  $B^i_{lm}$ and  $C^i_{lmn}$ let us introduce the following notations. Given whatever quantity ${\cal R}_{lmn}$ one can split it into {\it three} parts belonging to three representations of the cyclic group $Z_3$:
$$
{\cal R}^k_{lmn} ={ \cal R}^k_{(lmn)} +{\cal R}^k_{ \{lmn\} } +{\cal R}^k_{(lmn]} ,
$$
defined as follows:
\begin{eqnarray}
{\cal R}^k_{(lmn)}& =& \frac{1}{3} \, ( \, {\cal R}^k_{lmn} + {\cal R}^k_{nlm} + {\cal R}^k_{mnl} );\cr
{\cal R}^k_{ \{lmn \} }& =& \frac{1}{3} \, ( \, {\cal R}^k_{lmn} +  q^2 \, {\cal R}^k_{nlm} +
q \, {\cal R}^k_{mnl}) ;\cr
{\cal R}^k_{[lmn]} &=& \frac{1}{3} \, ( \, {\cal R}^k_{lmn} + q \, {\cal R}^k_{nlm} + q^2
\, {\cal R}^k_{mnl}) .
\end{eqnarray}
Now we can express $D ( d^2 x^k)$ in terms of coefficients $B^i_{lm}$
and $C^i_{lmn}$ as follows
\begin{equation}\label{expression through B and C}
D ( d^2 x^k)=B^i_{lm}\,Dx^l D^2 x^m +
    (C^i_{lms} - \Gamma^i_{r[s}\Gamma^r_{lm]}
    - \,\Gamma^r_{[ms}\Gamma^i_{l]r})
            \,Dx^l Dx^m Dx^s.
\end{equation}
Differentiating covariantly both sides of (\ref{trans of sec diff}) one
obtains the relation
\begin{equation}
D(d^2 x^k)=\frac{\partial x^k}{\partial y^{k'}}\,D(d^2 y^{k'})-
             \frac{\partial^2 x^k}{\partial y^{k'} \partial y^{l'}}
                \Gamma^{l'}_{r's'}Dx^{k'}Dx^{r'}Dx^{s'}.
\end{equation}
In order to give the transformation rules a more compact form
we shall use the following notations
$$
U^i_{j'}=\frac{\partial x^i}{\partial y^{j'}},\qquad
\partial_j  U^{i'}_{k}=\frac{\partial^2 y^{i'}}{\partial x^{k}\partial x^{j}}.
$$
Then replacing $D(d^2 x^k)$ and $D(d^2 y^{k'})$
in the above formula by their
expressions in terms of the coefficients $B^i_{lm}$ and $C^i_{lms}$
according to (\ref{expression through B and C}) and collecting together
similar terms we get the transformation rules
$$B^i_{lm} = B^{i'}_{l'm'}\; U^i_{i'}U^{l'}_{l} U^{m'}_{m}\  $$
$$C^i_{lms} = U^{i}_{i'} U^{l'}_{l}U^{m'}_{m}U^{n'}_{n}\, C^{i'}_{l'm's'}
+ U^{i}_{i'}U^{n'}_{[n}
 \partial_m    U^{r'}_{l]}\,\Gamma^{i'}_{r'n'} +
         U^i_{s'}\partial_ r U^{s'}_{[n} U^{l'}_{\underline{l}}U^{m'}_{m]}
U^{r}_{r'}\,\Gamma^{r'}_{l'm'}  $$
$$   + U^i_{s'}\partial _r U^{s'}_{[n}\partial _l  U^{t'}_{m]}U^{r}_{t'}+
 U^i_{i'}U^{l'}_{[n} \partial_l U^{r'}_{m]} \,\Gamma^{i'}_{l'r'}
     + U^{i}_{s'}\partial_r U^{s'}_{[n} U^{m'}_{l}U^{n'}_{m]}U^r_{r'}\,
         \Gamma^{r'}_{m'n'}+
         U^{i}_{s'}U^r_{t'}\partial_r U^{s'}_{[n}\partial_l U^{t'}_{m]} $$
As in the case of torsion in ordinary exterior calculus, the tensorial
character is displayed only by the irreducible part of these coefficients
displaying the corresponding symmetry.

Here is what we do mean by this. As in the usual case the connection
coefficients $\Gamma^k_{lm}$ could be split into two parts, the torsion
(antisymmetric part) and the symmetric part,
$$ \Gamma^k_{lm} = \frac{1}{2} \, [ \Gamma^k_{lm} + \Gamma^k_{ml} ] +
\frac{1}{2} \,[ \Gamma^k_{lm} -\Gamma^k_{ml}] =\Gamma^k_{(lm)} +S^k_{lm}.$$
so the four-index symbols $C^k_{lmn}$ as we have mentioned earlier can be split into {\it three} parts
belonging to three representations of the cyclic group $Z_3$ :
$$ C^k_{lmn} = C^k_{(lmn)} + C^k_{ \{lmn\} } + C^k_{(lmn]}. $$
Therefore, only the part ${\tilde{C}}^k_{[lmn]}$ has to be taken into account
in the final expression:
$$ D^3 x^k = {\tilde{B}}^k_{lm} \, D x^l D^2 x^m + {\tilde{C}}^k_{[lmn]} \,
D x^l D x^m D x^n = $$
\begin{equation}
= \biggl( B^k_{lm} + q^2 \, \Gamma^k_{ml} + q \, \Gamma^k_{lm} \biggr) \,
D x^l D^2 x^m + {\tilde{C}}^k_{[lmn]} \, D x^l D x^m D x^n
\end{equation}
with
$$ {\tilde{C}}^k_{lmn} = C^k_{lmn} + \partial_l \Gamma^k_{mn} - \Gamma^r_{lm}
\Gamma^k_{rn} - q \, \Gamma^r_{mn} \Gamma^k_{lr} $$
Because the coefficients in both terms on the right-hand side are tensors,
we can start to investigate their intrinsic properties. Among these, the
condition of reality should be applied to both coefficients separately.
Starting with the first coefficient, ${\tilde{B}}^k_{lm}$, and recalling that
\vskip 0.2cm
\centerline{$ q = - \frac{1}{2} + \frac{i \sqrt{3}}{2} \, \ \ \, {\rm and} \, \ \ \,
 q^2 = - \frac{1}{2} - \frac{i \sqrt{3}}{2} $}
\vskip 0.2cm
\noindent
we have explicitly
$${\tilde{B}}^k_{lm} = B^k_{lm} - \frac{1}{2} \, \biggl( \Gamma^k_{ml} +
\Gamma^k_{lm} \biggr) \, + \frac{i \sqrt{3}}{2} \, \biggl( \Gamma^k_{lm}
- \Gamma^k_{ml} \, \biggr) $$
The imaginary part is the {\it torsion tensor} of the connection
$\Gamma^k_{lm}$, so the reality of the coefficient ${\tilde{B}}^k_{lm}$ is
equivalent with the vanishing of the torsion, leaving only the {\it symmetric}
part of $\Gamma^k_{lm}$. So, from now on, we can write
$${\tilde{B}}^k_{lm} =  B^k_{lm} - \Gamma^k_{lm}, \, \ \ \, {\rm with} \, \ \
\, \Gamma^k_{lm} = \Gamma^k_{ml}. $$
This means that the $B^k_{lm}$ transform as connection coefficients, so that
the difference $B^k_{lm} - \Gamma^k_{lm}$ is a tensor. As a corollary, the
vanishing of $D^3 x^k$ implies that $B^k_{lm} = \Gamma^k_{lm}$ and
$\Gamma^k_{lm} = \Gamma^k_{ml}$.

The symmetry of the connection coefficients makes it possible to identify
the tensor appearing in the coefficients ${\tilde{C}}^k_{lmn}$. As a matter
of fact, only the part ${\tilde{C}}^k_{[lmn]}$ is relevant here, the other
two irreducible parts' contribution vanishing when contracted with the
covariant expression $D x^l  D x^m D x^n$. The part of ${\tilde{C}}^k_{lmn}$
containing the connection coefficients and their derivatives should be also
$Z_3$-anti-symmetrized, yielding
$$ \frac{1}{3} \, \biggl( C^k_{lmn} + q C^k_{nlm} + q^2 C^k_{mnl} +
\partial_l \Gamma^k_{mn} + q \partial_n \Gamma^k_{lm} +  q^2 \partial_m
\biggr) \, D x^l D x^m D x^n
\Gamma^k_{nl}  $$
$$ - \frac{1}{3} \, \biggl( \Gamma^r_{lm} \Gamma^k_{rn} - q \Gamma^r_{nl}
\Gamma^k_{rm} - q^2 \Gamma^r_{mn} \Gamma^k_{rl} - q \Gamma^r_{mn} - q^2
\Gamma^r_{lm} \Gamma^k_{nr} - \Gamma^r_{nl} \Gamma^k_{rm} \biggr) \,
D x^l D x^m D x^n $$
\indent
It is not difficult, taking into account the symmetries, to identify the final
result in terms of the Riemann tensor:
${\tilde{C}}^k_{lmn} \, D x^l D x^m D x^n $ is equal to
$$ \biggl( C^k_{[lmn]} + \frac{1}{3} \, [ R^k_{ \, nlm} + R^k_{\, mln} \, ]
+ \frac{q}{3} \, [ R^k_{ \, mnl} + R^k_{\, lnm} \, ] + \frac{q^2}{3} \,
[ R^k_{\, lmn} + R^k_{nml} \, ]  \, \biggr) \, D x^l D x^m D x^n $$
\indent
Taking into account that
$$C^k_{[lmn]} = \frac{1}{3} \, \biggl( C^k_{[lmn]} + q \, C^k_{[nlm]} +
q^2 \, C^k_{[mnl]} \biggr) $$
and assuming that the coefficients $C^k_{[lmn]}$ are real, the vanishing of
the above expression leads to the equality
\begin{equation}
C^k_{[lmn]} = - \, [ \, R^k_{\, nlm} \, + \, R^k_{\, mln} \, ]
\end{equation}
\indent
The analogy with the usual exterior differential calculus is now obvious. In
the usual case the condition $D^2 x^k = 0$ implied the vanishing of torsion,
$S^k_{lm} = \frac{1}{2} [\Gamma^k_{lm} - \Gamma^k_{ml}] = 0$, whereas now, in
our $Z_3$-graded case, the similar condition $D^3 x^k = 0$ implies not only
the vanishing of torsion, but also determines entirely the coefficients
$B^k_{lm}$ (equal to  $\Gamma^k_{(ml)}$), and partly the coefficients
$C^k_{lmn}$ namely, their $q$-skew-symmetric part $C^k_{[lmn]}$ (equal then
to the expression $ - [R^k_{\,nlm} + R^k_{\, mln}]$). By analogy, one can
impose similar conditions on the ``conjugate'' $q^2$-skewsymmetric part
$C^k_{ \{ mnl \} }$ , defining it e.g. as $C^k_{ \{ mnl \} } = C^k_{[lnm]}$.
However, the {\it totally symmetric} part $C^k_{(lmn)}$, which is not a
tensor, is still undefined, because its contribution cancels automatically
when contracted with the $3$-form $ D x^l D x^m D x^n$.

The symmetric part of $C^k_{lmn}$ together with the coefficients $B^k_{lm}$
may be used for the definition of a new kind of parallel transport and
generalized geodesic curves. One can define, independently of usual covariant
derivative of a vector along a parametrized curve $x^k(\lambda)$ determined
by given connection coefficients $\Gamma^k_{lm}$,
\begin{equation}
\frac{D Y^k}{D \lambda} = \frac{d Y^k}{d \lambda} + \Gamma^k_{lm} \,
\frac{d x^l}{d \lambda} Y^m = 0 ,
\end{equation}
a second-order covariant derivative which is not an iteration of the first one:
\begin{equation}
\frac{{\cal{D}}^2 Y^k}{{\cal{D}} \lambda^2} = \frac{d^2 Y^k}{d \lambda^2} +
E^k_{lm} \, \frac{d x^k}{d \lambda} \frac{D Y^m}{D \lambda} + F^k_{lm} \,
\frac{D^2 x^l}{D \lambda^2} \, Y^m + G^k_{lmn} \, \frac{d x^l}{d \lambda}
\frac{d x^m}{d \lambda} \, Y^n
\end{equation}
where we use a different notation, $\frac{{\cal{D}} \,}{{\cal{D}} \lambda}$
in order to stress that we don't consider here a simple iteration of the usual
covariant differentiation $\frac{D \, }{D \lambda}$. The coefficients
$E^k_{lm}$ and  $F^k_{lm}$ are not identical a priori; all we need to know
about the transformation properties of $E^k_{lm}, \, F^k_{lm}$ and $G^k_{lmn}$
is that the resulting quantity $\frac{{\cal{D}}^2 Y^k}{{\cal{D}} \lambda^2}$
transforms as a vector under a coordinate change.

If we replace the vector field $Y^k \, (x^m(\lambda))$ by the vector
$\frac{d x^k}{d \lambda}$ tangent to the curve, we obtain a third-order
generalization of the geodesic equation:
\begin{equation}
\frac{{\cal{D}}^3 x^k}{{\cal{D}} \lambda^3} = \frac{d^3 x^k}{d \lambda^3} +
[ \,E^k_{lm} + F^k_{ml} \,] \, \frac{d x^l}{d \lambda}
\frac{D^2 x^m}{D \lambda^2} + G^k_{lmn} \, \frac{d x^l}{d \lambda}
\frac{d x^m}{d \lambda} \frac{d x^n}{d \lambda} = 0
\label{3geodesic}
\end{equation}
\indent
Now the $\frac{d x^k}{d \lambda}$ etc. are commutative entities, so that
$G^k_{lmn} = C^k_{(lmn)}$. Had we iterated the usual covariant derivative
in the above equation, then the coefficients $E^k_{lm} + F^k_{lm}$ and
$C^k_{(lmn)}$ would be completely determined from the connection coefficients
$\Gamma^k_{lm}$ and their derivatives; however, we can introduce more general
coefficients having the required transformation properties and independent of
$\Gamma^k_{lm}$. This is reminiscent of a similar situation one level below,
when the Christoffel connection is totally determined by the metric, but a
larger class of affine connections exist which are independent of metric.

The generalized geodesic equation of third order (\ref{3geodesic}) defines
a larger class of curves that the usual geodesics and may be of interest in
probing certain geometrical objects. For example, in the flat Euclidean space
the solutions of (\ref{3geodesic}) include not only the straight lines, but
also all possible hyperbolae.

The geometric aspects of the new differential calculus are beyond the scope
of our present article, and we shall publish them later.
\newline
\vskip .4cm \noindent {\tbf ACKNOWLEDGMENTS} \vskip .3cm The
authors are grateful to M. Dubois-Violette, O. Suzuki and L.
Vainerman for valuable discussions. The first author (V.A.) would
like to acknowledge the financial support of the Estonian Science
Foundation under the grant No. 2403.

\end{document}